\documentclass[reqno]{amsart}
\usepackage{amssymb}
\usepackage{amsmath}
\usepackage{amsfonts}
\usepackage{amscd}
\usepackage[small,nohug,heads=vee]{diagrams}

\diagramstyle[labelstyle=\scriptstyle]
\newtheorem{theorem}{Theorem}
\theoremstyle{plain}

\newtheorem{definition}[theorem]{Definition}

\newtheorem{problem}[theorem]{Problem}

\newtheorem{remark}[theorem]{Remark}

\numberwithin{equation}{section}
\numberwithin{theorem}{subsection}
\input{tcilatex}

\begin{document}
\title[The categorical Weil representation]{The categorical Weil
representation}
\author{Shamgar Gurevich}
\curraddr{Department of Mathematics, University of Wisconsin, Madison, WI
53706, USA.}
\email{shamgar@math.wisc.edu}
\author{Ronny Hadani}
\curraddr{Department of Mathematics, The University of Texas at Austin,
Austin, TX 78712, USA. }
\email{hadani@math.utexas.edu}
\thanks{\copyright\ Copyright by S. Gurevich and R. Hadani, November 2010.
All rights reserved.}

\begin{abstract}
In a previous work the authors gave a conceptual explanation for the
linearity of the Weil representation over a finite field $k$ of odd
characteristic: There exists a canonical system of intertwining operators
between the Lagrangian models of the Heisenberg representation. This defines
a canonical vector space $\mathcal{H(}V\mathcal{)}$ associated with a
symplectic vector space $V$ over $k.$ In this paper we prove a general
theorem about idempotents in categories, and we use it to solve the sign
problem, formulated by Bernstein and Deligne, on the compatibility between
the associativity constraint and the convolution structure of the $\ell $%
-adic sheaf of canonical intertwining kernels. This sheaf governs---via the
sheaf-to-function correspondence---the function theoretic system of
intertwiners. As an application we define a canonical category $\mathcal{C}(%
\mathbf{V})$ associated with the symplectic vector space variety $\mathbf{V}$%
, and we obtain the canonical model of the categorical Weil representation.
\end{abstract}

\maketitle

\section{Introduction}

\subsection{The Weil representation}

The Weil representation \cite{W} over a finite field $k$ is the algebra
object that governs the symmetries of the standard Hilbert space $\mathcal{H}%
_{L}=L^{2}(L,%
\mathbb{C}
)$ of complex valued functions on a finite-dimensional vector space $L$ over 
$k.$ In this paper we will be interested only with the case where $k$ is of
odd characteristic. We have the (split) symplectic vector space $V=L\times
L^{\ast }$ with its standard symplectic form. We denote by $Sp(V)$ the
corresponding symplectic group. A nontrivial argument, due to Schur, implies
that there exists a linear representation\bigskip 
\begin{equation*}
\rho _{L}:Sp(V)\rightarrow GL(\mathcal{H}_{L}),
\end{equation*}%
obtaining the Weil representation. This representation depends also on a
choice of additive character $\psi :k\rightarrow 
\mathbb{C}
^{\ast }.$

The Weil representation is a central object of modern harmonic analysis and
the theory of the discrete Fourier transform. It has many applications in
automorphic forms, number theory, the theory of theta functions,
mathematical physics, coding theory, signal processing, and other domains of
knowledge.

Probably, the most interesting fact, for researchers who are not familiar
with the Weil representation, is that it can be thought of as a group of
operators which includes the discrete Fourier transform. Indeed, for any
nondegenerate symmetric bilinear form $B$ on $L$ we have \cite{D1}%
\begin{equation*}
\left[ \rho _{L}\left( 
\begin{pmatrix}
0 & -B^{-1} \\ 
B & \text{ \ }0%
\end{pmatrix}%
\right) f\right] (x)=\frac{1}{G(B,\psi )}\sum_{y\in L}\psi (B(x,y))f(y),
\end{equation*}%
where $G(B,\psi )$ is an appropriate Gauss sum normalization.

\subsection{Canonical vector space}

In \cite{GH1, GH2} a conceptual explanation for the existence of the linear
Weil representation was proposed. Specifically, it was shown that there
exists an explicit quantization functor 
\begin{equation}
\mathcal{H}:\mathsf{Symp}\rightarrow \mathsf{Vect},  \label{Q}
\end{equation}%
where $\mathsf{Symp}$ denotes the (groupoid) category whose objects are
finite dimensional symplectic vector spaces over the finite field $k$, and
morphisms are linear isomorphisms of symplectic vector spaces and $\mathsf{%
Vect}$ denotes the category of finite dimensional complex vector spaces.

As a consequence, for a fixed symplectic vector space $V\in $ $\mathsf{Symp}$%
, we obtain, by functoriality, a homomorphism 
\begin{equation}
\rho _{V}:Sp(V)\rightarrow GL(\mathcal{H}(V)),  \label{CWeil}
\end{equation}%
which is isomorphic to the Weil representation. We refer to the vector space 
$\mathcal{H}(V)$ as the \textit{canonical vector space }associated to $V$,
and to the representation (\ref{CWeil}) as the \textit{canonical model} of
the Weil representation of $Sp(V)$.

\subsection{Main results}

\subsubsection{\textbf{Canonical category}}

In this paper, we obtain the categorical analog of the functor (\ref{Q}),
i.e., we define an explicit quantization lax $2$-functor \cite{Gr1, V} 
\begin{equation*}
\mathbf{V\mapsto }\text{ }\mathcal{C}(\mathbf{V}),
\end{equation*}%
associating certain category of $\ell $-adic sheaves $\mathcal{C}(\mathbf{V})
$ to any object $\mathbf{V}$ in the (groupoid) category $\mathbf{Symp}$
whose objects are finite dimensional symplectic vector spaces in the
category of algebraic varieties over $k$, and morphisms are linear
isomorphisms of symplectic vector spaces. We will refer to the category $%
\mathcal{C}(\mathbf{V})$ as the \textit{canonical category }associated to $%
\mathbf{V.}$

\subsubsection{\textbf{The categorical Weil representation}}

In particular, for a fixed object $\mathbf{V\in Symp}$ we obtain, by
functoriality, an action of the algebraic group $\mathbf{Sp=}Sp(\mathbf{V})$
on the category $\mathcal{C}(\mathbf{V}),$ forming the categorical analog of
(\ref{CWeil}). For the sake of the introduction, it is enough to say that
this action, which we will call the \textit{categorical Weil representation}%
, is an explicit \textit{family} of auto-equivalence functors 
\begin{equation*}
\mathbf{\rho }_{\mathbf{V}}(g)\mathbf{:\mathcal{C}(\mathbf{V})\rightarrow 
\mathcal{C(}V)},\text{ \ }g\in \mathbf{Sp,}
\end{equation*}%
which are induced from the action morphism $\mathbf{Sp\times V\rightarrow V.}
$

\subsubsection{\textbf{The idempotent theorem}}

The main technical result of this paper is the proof of the idempotent
theorem which is needed in order to define the category $\mathbf{\mathcal{C(}%
V).}$ This theorem proposes a solution to a problem, formulated by Bernstein 
\cite{B} and Deligne \cite{D2}, that we will call the \textit{sign problem}.
We devote the rest of the introduction to an intuitive description of the
sign problem, its solution, and its implication to the definition of $%
\mathbf{\mathcal{C(}V).}$

\subsection{Canonical system of intertwining operators}

The existence of the canonical vector space $\mathcal{H}(V)$ is a
manifestation of the existence of a canonical system of intertwining
operators, between the models of the Heisenberg representation, associated
with oriented Lagrangian subspaces in $V$.

We denote by $OLag=OLag(V)$ the set of oriented Lagrangian subspaces in $V$,
i.e., the set of pairs $L^{\circ }=(L,o_{L})$ where $L\subset V$ is a
Lagrangian subspace and $o_{L}$ is a nonzero vector in the top wedge product 
$\tbigwedge^{top}L.$ In addition, we denote by $H=H(V)$ the Heisenberg group
associated with $V$, and by $Z=Z(H)\simeq k$ its center. Finally, let us
choose a nontrivial additive character $\psi :Z\rightarrow 
\mathbb{C}
^{\ast }.$

The Stone--von Neumann theorem asserts that there exist a unique (up to
isomorphism) irreducible representation $\pi :H\rightarrow GL(\mathcal{H})$
with $\pi (z)=\psi (z)\cdot Id_{\mathcal{H}}$ for every $z\in Z.$ We will
call the representation $(\pi ,H,\mathcal{H)}$ the \textit{Heisenberg
representation}. An important family of models of the Heisenberg
representation is associated with oriented Lagrangian subspace in $V$. To
every $L^{\circ }\in OLag$ we have the model $(\pi _{L^{\circ }},H,\mathcal{H%
}_{L^{\circ }})$, where the vector space $\mathcal{H}_{L^{\circ }}$ is the
space $%
\mathbb{C}
\left( L\backslash H,\psi \right) $ of functions $f:H\rightarrow 
\mathbb{C}
$ such that $f\left( z\cdot l\cdot h\right) =\psi \left( z\right) f\left(
h\right) $, for every $z\in Z$, $l\in L,$ and the action $\pi _{L^{\circ }}$
is given by right translation. The collection of models $\left\{ \mathcal{H}%
_{L^{\circ }}\right\} $ can be thought of as a vector bundle $\mathfrak{H}$
on the set $OLag$, with fibers $\mathfrak{H}_{|L^{\circ }}=\mathcal{H}%
_{L^{\circ }}$. \ The main technical result proved in \cite{GH1, GH2} is the
strong Stone--von Neumann theorem, i.e., the vector bundle $\mathfrak{H}$
admits a canonical trivialization. Concretely, this means that there exists
a canonical system of intertwining operators $T_{M^{\circ },L^{\circ }}\in 
\mathrm{Hom}_{H}\left( \mathcal{H}_{L^{\circ }},\mathcal{H}_{M^{\circ
}}\right) $, for every $M^{\circ },L^{\circ }\in OLag$, satisfying the
following multiplicativity property: 
\begin{equation}
T_{N^{\circ },M^{\circ }}\circ T_{M^{\circ },L^{\circ }}=T_{N^{\circ
},L^{\circ }},  \label{Mult-Int}
\end{equation}%
for every $N^{\circ },M^{\circ },L^{\circ }\in OLag$.

\begin{remark}
Similar considerations with $k=%
\mathbb{R}
$ yield \cite{LV} a canonical system of intertwiners, which are
multiplicative up to a sign $\pm 1,$ referred to as the metaplectic sign.
\end{remark}

The canonical vector space $\mathcal{H}\left( V\right) $ is the space of
"horizontal sections" of $\mathfrak{H}$%
\begin{equation*}
\mathcal{H}\left( V\right) =\Gamma _{hor}\left( OLag,\mathfrak{H}\right) ,
\end{equation*}%
i.e., a vector in $\mathcal{H}\left( V\right) $ is a compatible system $%
\left( f_{L^{\circ }}\in \mathcal{H}_{L^{\circ }};\text{ }L^{\circ }\in
OLag\right) $ such that $T_{M^{\circ },L^{\circ }}\left( f_{L^{\circ
}}\right) =f_{M^{\circ }}$, for every $M^{\circ },L^{\circ }\in OLag$. The
symplectic group $Sp\left( V\right) $ acts on the vector space $\mathcal{H}%
\left( V\right) $ in an obvious manner and we obtain the model (\ref{CWeil}).

In order to formulate the sign problem \cite{B, D2}, we will need the
description of the system $\{T_{M^{\circ },L^{\circ }}\}$ by kernels and
their geometrization.

\subsection{Canonical system of intertwining kernels}

Every intertwining operator $T_{M^{\circ },L^{\circ }}$ can be uniquely
presented by a kernel function $K_{M^{\circ },L^{\circ }}\in 
\mathbb{C}
\left( M\backslash H/L,\psi \right) $. The collection of kernel functions $%
\left\{ K_{M^{\circ },L^{\circ }}\right\} $ can be thought of as a single
function%
\begin{equation}
K\in 
\mathbb{C}
(OLag^{2}\times H),  \label{Ker-Fun}
\end{equation}%
given by $K\left( M^{\circ },L^{\circ },-\right) =K_{M^{\circ },L^{\circ
}}\left( -\right) $, for every $\left( M^{\circ },L^{\circ }\right) \in
OLag^{2}$. \ The multiplicativity property (\ref{Mult-Int}) transformed into
the following convolution property with respect to the convolution $\ast $
of functions on the Heisenberg group:%
\begin{equation}
p_{32}^{\ast }K\ast p_{21}^{\ast }K=p_{31}^{\ast }K,  \label{Conv-Ker}
\end{equation}%
where $p_{ji}:OLag^{3}\times H\rightarrow OLag^{2}\times H$ are the
projections given by $p_{ji}(L_{3}^{\circ },L_{2}^{\circ },L_{1}^{\circ
},h)=(L_{j}^{\circ },L_{i}^{\circ },h)$, for $1\leq i<j\leq 3.$

We proceeds to the geometrization of the kernels.

\subsection{The sheaf of canonical geometric intertwining kernels}

\subsubsection{Geometrization}

A general ideology due to Grothendieck is that any meaningful set-theoretic
object is governed by a more fundamental algebra-geometric one. \textit{\ }%
The procedure by which one translate from the set theoretic setting to
algebraic geometry is called \textit{geometrization}, which is a formal
procedure by which sets are replaced by algebraic varieties and functions
are replaced by certain sheaf-theoretic objects.

The precise setting consists of:

\begin{itemize}
\item A set $X=\mathbf{X}(k)$ of rational points of an algebraic variety $%
\mathbf{X}$ defined over $k$.

\item A complex valued function $f\in 
\mathbb{C}
\left( X\right) $ governed by an $\ell $-adic Weil sheaf $\mathcal{F}$ on $%
\mathbf{X}$.
\end{itemize}

The variety $\mathbf{X}$ is a space equipped with an endomorphism $Fr:%
\mathbf{X}\rightarrow \mathbf{X}$, called Frobenius, such that the set $X$
is naturally identified with the set of fixed points $X=\mathbf{X}^{Fr}$. \ 

The sheaf $\mathcal{F}$ can be thought of as a vector bundle on the variety $%
\mathbf{X}$, equipped with an endomorphism $\vartheta :\mathcal{F\rightarrow
F}$ which lifts $Fr$.

The relation between the function $f$ and the sheaf $\mathcal{F}$ is called
Grothendieck's \textit{sheaf-to-function correspondence: }Given a point $%
x\in X$, the endomorphism $\vartheta $ restricts to an endomorphism $%
\vartheta _{x}:\mathcal{F}_{|x}\rightarrow \mathcal{F}_{|x}$ of the fiber $%
\mathcal{F}_{|x}$. The value of $f$ on the point $x$ is given by 
\begin{equation*}
f(x)=f^{\mathcal{F}}\left( x\right) =Tr(\vartheta _{x}:\mathcal{F}%
_{|x}\rightarrow \mathcal{F}_{|x}).
\end{equation*}

\subsubsection{The sheaf of canonical geometric intertwining kernels}

The function $K$ (\ref{Ker-Fun}) fits nicely to the geometrization
procedure. Denote by $\mathbf{OLag}^{2}\mathbf{\times H}$ the algebraic
variety with $OLag^{2}\times H=\mathbf{OLag}^{2}\mathbf{\times H(}k\mathbf{).%
}$ In \cite{GH1} the authors defined a geometrically irreducible, (shifted)
perverse, $\ell $-adic Weil sheaf $\mathcal{K}$ on $\mathbf{OLag}^{2}\mathbf{%
\times H,}$ that satisfies the following two properties:

\begin{enumerate}
\item \textbf{Convolution. }There exists a canonical isomorphism $\theta :%
\mathcal{K\ast K}\widetilde{\mathcal{\rightarrow }}\mathcal{K}$.\textbf{\ }
\ 

\item \textbf{Function. }Applying sheaf-to-function procedure we recover $f^{%
\mathcal{K}}=K.$
\end{enumerate}

Here, the notation $\theta :\mathcal{K\ast K}\widetilde{\mathcal{\rightarrow 
}}\mathcal{K}$ stands to simplify the more precise notation which is the
geometric analogue of (\ref{Conv-Ker}). \ We will call $\mathcal{K}$ the 
\textit{sheaf of canonical geometric intertwining kernels.}

We are ready now to formulate the sign problem.

\subsection{The sign problem}

In an attempt to understand the nature of the ("disappearance" of the)
metaplectic sign over finite fields of odd characteristic, we formulate \cite%
{B, D2} the following sign problem. Let $\mathcal{K}$ be the sheaf of
canonical geometric intertwining kernels. Consider the commutative diagram

\begin{equation}
\begin{CD} (\mathcal{K} \ast \mathcal{K}) \ast \mathcal{K} @>{\alpha}>>
\mathcal{K} \ast (\mathcal{K} \ast \mathcal{K})\\ @VV{\theta \, \ast \, id}
V @VV{id \, \ast \, \theta}V\\ \mathcal{K} \ast \mathcal{K} @. \mathcal{K}
\ast \mathcal{K}\\ @VV{\theta}V @VV{\theta}V\\ \mathcal{K} @>{C}>>
\mathcal{K} \end{CD}  \label{BP}
\end{equation}%
where $\alpha $ is the \textit{associativity constraint}\ \cite{DMOS} for
convolution $\ast $ of sheaves on the Heisenberg group $\mathbf{H}$, the
morphism $\theta $ is the isomorphism appearing in the convolution property
of the sheaf $\mathcal{K}$, and $C$ is by definition the isomorphism that
makes the diagram commutative.

The sheaf $\mathcal{K}$ is geometrically irreducible, hence, $C=c\cdot id$
is a scalar morphism.\smallskip

\begin{problem}[The sign problem]
Compute the scalar $c.\smallskip $
\end{problem}

The idea is \cite{B, D2} that the value of $c$ might suggests a new
understanding about the nature of the metaplectic sign over finite field,
i.e., that it "moves" one hierarchy higher, becoming a property of the sheaf 
$\mathcal{K}$ that cannot be observed on the level of the function $K$.

\subsection{A solution to the sign problem}

The main technical result of this paper is the proof of the following
theorem:

\begin{theorem}[The idempotent theorem---particular case]
\label{IT}We have $c=1.$
\end{theorem}

In fact, we will show that the \textit{idempotent theorem }holds in a very
general situation. Let $\mathcal{C}$ be any category with an "operation",
i.e., a functor $\otimes :\mathcal{C\times C\rightarrow C}$ equipped with
associativity constraint $\alpha $ \cite{DMOS}, and suppose $P\in \mathcal{C}
$ is an "idempotent", i.e., equipped with an isomorphism $\theta :P\otimes P%
\widetilde{\rightarrow }P$. Then we can form an analogue diagram to (\ref{BP}%
) and obtain an isomorphism $C:P\rightarrow P.$ We will prove that under
natural conditions $C=id.$

Let us describe shortly the application of Theorem \ref{IT} to the
definition of the canonical category.

\subsection{ The canonical category}

The idempotent theorem suggests the definition of a canonical category $%
\mathcal{C(}\mathbf{V}\mathcal{)}$ associated with the symplectic vector
space variety $\mathbf{V}$ with $V=\mathbf{V(}k\mathbf{).}$ This category is
the categorical analogue of the vector space $\mathcal{H}(V),$ and it is the
basic object behind the existence of the \textit{categorical Weil
representation}. The category $\mathcal{C(}\mathbf{V}\mathcal{)}$ consists
of $\ell $-adic sheaves on $\mathbf{OLag\times H}$ satisfying geometric
conditions that are analog of these satisfied by vectors in $\mathcal{H}(V)$%
. The most important condition is that each $\mathcal{F\in C(}\mathbf{V}%
\mathcal{)}$ is equipped with an isomorphism 
\begin{equation*}
\eta :\mathcal{K\ast }p_{1}^{\ast }\mathcal{F}\widetilde{\rightarrow }%
p_{2}^{\ast }\mathcal{F}_{,}
\end{equation*}%
which is compatible with $\alpha $ and $\theta $. Here, $p_{i}:\mathbf{OLag}%
^{2}\mathbf{\times H\rightarrow OLag\times H}$ are given by $%
p_{i}(L_{2}^{\circ },L_{1}^{\circ },h)=(L_{i}^{\circ },h),$ $i=1,2,$ and
compatibility means the commutativity 
\begin{equation*}
\begin{CD} (\mathcal{K} \ast \mathcal{K}) \ast \mathcal{F} @>{\alpha}>>
\mathcal{K} \ast (\mathcal{K} \ast \mathcal{F})\\ @VV{\theta \, \ast \, id}
V @VV{id \, \ast \, \eta}V\\ \mathcal{K} \ast \mathcal{F} @. \mathcal{K}
\ast \mathcal{F}\\ @VV{\eta}V @VV{\eta}V\\ \mathcal{F} @>{id}>> \mathcal{F}
\end{CD}
\end{equation*}%
This definition make sense only if $C=id$ in (\ref{BP}).

\subsection{Structure of the paper}

In Section \ref{CVS} we recall the construction of the canonical vector
space, and the Weil representation, using the canonical intertwining
operators which we describe in the language of kernels. In Section \ref{GCIK}
we recall the sheaf theoretic counterpart of the intertwining kernels. In
Section \ref{TIT} we define the notion of idempotent in a category, and we
state the \textit{idempotent theorem} (Theorem \ref{ITG}) on the
compatibility between an idempotent structure and an associativity
constraint for operations in categories. We describe also a generalization
of the idempotent theorem which suggests a solution to the sign problem. In
Section \ref{Appl} we describe the following applications: We formulate and
solve the sign problem; we define the canonical category which is a
categorification of the canonical vector space; we obtain the canonical
model of the categorical Weil representation; and we suggest a solution to
the sign problem in the context of the geometric Weil representation sheaf.
Finally, in Appendix \ref{PITG} we suggest a proof for the idempotent
theorem.\medskip

\textbf{Acknowledgement. }We thank our teacher J. Bernstein who initiated
and guided the research reported in this paper. We thank O. Gabber for his
important contributions, and for his invitation to visit IHES, Paris, in
2009, and 2011, to discuss parts of the content of this paper. Finally, we
acknowledge P. Deligne for several communications, discussions, and meetings
at IAS, Princeton, during 2008--2011.

\section{The canonical vector space\label{CVS}}

We would like to recall the construction of the canonical vector space
established in \cite{GH1, GH2}, to describe the family of canonical
intertwining operators and their explicit presentation as kernels.

\subsection{The Heisenberg representation \label{HR}}

Let $(V,\omega )$ be a $2n$--dimensional symplectic vector space over the
finite field $k=\mathbb{F}_{q}$, where $q$ is odd. The \textit{Heisenberg }%
group associated with $V$ can be presented as the set $H=V\times k$ with the
multiplication given by%
\begin{equation*}
(v,z)\cdot (v^{\prime },z^{\prime })=(v+v^{\prime },z+z^{\prime }+\tfrac{1}{2%
}\omega (v,v^{\prime })).
\end{equation*}

The center of $H$ is $\ Z=Z(H)=\{(0,z):$ $z\in k\}.$

Let $\psi :Z\rightarrow 
\mathbb{C}
^{\ast }$ be a nontrivial character of the center.

\begin{theorem}[Stone--von Neuman]
\label{S-vN_thm}There exists a unique (up to isomorphism) irreducible
unitary representation $(\pi ,H,\mathcal{H)}$ with the center acting by $%
\psi ,$ i.e., $\pi _{|Z}=\psi \cdot Id_{\mathcal{H}}$.
\end{theorem}

The representation $\pi $ which appears in the above theorem will be called
the \textit{Heisenberg representation}.

\subsection{Oriented Lagrangian models}

We are interested in a particular family of models of the Heisenberg
representation which are associated with oriented Lagrangian subspaces in $V$%
. Let $Lag$ denotes the set of Lagrangian subspaces in $V.$

\begin{definition}
An \underline{oriented Lagrangian subspace} in $V$ is a pair $L^{\circ
}=\left( L,o_{L}\right) $, where $L\in Lag$ and $o_{L}\in $ $%
\bigwedge\nolimits^{top}L$ is a nonzero vector.
\end{definition}

Let us denote by $OLag$ the set of oriented Lagrangian subspaces in $V.$ We
associate with each oriented Lagrangian $L^{\circ }\in OLag$ a model $(\pi
_{L^{\circ }},H,\mathcal{H}_{L^{\circ }})$ of the Heisenberg representation
as follows: The vector space $\mathcal{H}_{L^{\circ }}$ is the space $%
\mathbb{C}
(L\backslash H,\psi )$ of functions $f:H\rightarrow 
\mathbb{C}
$ satisfying $f\left( z\cdot l\cdot h\right) =\psi (z)f\left( h\right) $ for
every $z\in Z,$ $l\in L,$ and the Heisenberg action is given by right
translation $\pi _{L^{\circ }}\left( h\right) [f]\left( h^{\prime }\right)
=f\left( h^{\prime }\cdot h\right) $ for every $f\in \mathcal{H}_{L^{\circ
}} $.

\subsection{System of canonical intertwining operators \label{SS-vN_sub}}

The collection of models $\left\{ \mathcal{H}_{L^{\circ }}\right\} $ forms a
vector bundle $\mathfrak{H}\rightarrow OLag$ with fibers $\mathfrak{H}%
_{L^{\circ }}=\mathcal{H}_{L^{\circ }}$ acted upon by $H$ via $\pi
_{L^{\circ }}.$

\begin{definition}
Let $\mathfrak{E}\rightarrow OLag$ be an $H$-vector bundle. A \underline{%
trivialization} of $\mathfrak{E}$ is a system of intertwining operators
(isomorphisms) $\{E_{M^{\circ },L^{\circ }}\in \mathrm{Hom}_{H}(\mathfrak{E}%
_{L^{\circ }},\mathfrak{E}_{M^{\circ }}):\left( M^{\circ },L^{\circ }\right)
\in OLag^{2}\}$ satisfying the following multiplicativity condition%
\begin{equation*}
E_{N^{\circ },M^{\circ }}\circ E_{M^{\circ },L^{\circ }}=E_{N^{\circ
},L^{\circ }},
\end{equation*}%
for every $N^{\circ },M^{\circ },L^{\circ }\in OLag$.
\end{definition}

The main result of \cite{GH1, GH2} is the following:

\begin{theorem}[The strong S-vN property]
\label{SS-vN_thm}The $H$-vector bundle $\mathfrak{H}$ admits a natural
trivialization $\{T_{M^{\circ },L^{\circ }}\}$.
\end{theorem}

The intertwining operators $\{T_{M^{\circ },L^{\circ }}\}$ in the above
theorem will be referred to as the system of \textit{canonical intertwining
operators}.

\subsection{Explicit formulas for the canonical intertwining operators}

Let us denote by $U_{2}\subset OLag\left( V\right) ^{2}$ the subset
consisting of pairs of oriented Lagrangians $\left( M^{\circ },L^{\circ
}\right) \in OLag\left( V\right) ^{2}$ which are in general position, i.e., $%
L+M=V$. For every $\left( M^{\circ },L^{\circ }\right) \in U_{2}$, we have 
\cite{GH1} 
\begin{equation*}
T_{M^{\circ },L^{\circ }}=A_{M^{\circ },L^{\circ }}\cdot F_{M^{\circ
},L^{\circ }},
\end{equation*}%
where $F_{M^{\circ },L^{\circ }}:\mathcal{H}_{L^{\circ }}\rightarrow 
\mathcal{H}_{M^{\circ }}$ is the averaging morphism given by%
\begin{equation*}
F_{M^{\circ },L^{\circ }}\left[ f\right] \left( h\right) =\sum\limits_{m\in
M}f\left( m\cdot h\right) ,
\end{equation*}%
for every $f\in \mathcal{H}_{L^{\circ }}$ and $A_{M^{\circ },L^{\circ }}$ is
a normalization constant given by 
\begin{equation}
A_{M^{\circ },L^{\circ }}=\left( G(\psi )/q\right) ^{n}\sigma (\left(
-1\right) ^{\left( \QATOP{n}{2}\right) }\omega _{\wedge }\left(
o_{L},o_{M}\right) ),  \label{A}
\end{equation}

where

\begin{itemize}
\item $\sigma $ is the unique quadratic character (also called the Legendre
character) of the multiplicative group $G_{m}=k^{\ast }$.

\item $G(\psi )$ is the one-dimensional Gauss sum 
\begin{equation*}
G(\psi )=\sum\limits_{z\in k}\psi (\tfrac{1}{2}z^{2}).
\end{equation*}

\item $\omega _{\wedge }:\bigwedge\nolimits^{top}L\times $ $%
\bigwedge\nolimits^{top}M\rightarrow $ $k$ is the pairing induced by the
symplectic form.
\end{itemize}

\subsection{The canonical vector space and the Weil representation\label%
{CH_sub}}

Using Theorem \ref{SS-vN_thm}, we can associate, in a functorial manner, a
vector space $\mathcal{H}\left( V\right) $ \ to the symplectic vector space $%
V$ as follows: Define $\mathcal{H}\left( V\right) $ to be the space of
"horizontal sections" of the trivialized vector bundle $\mathfrak{H}$ 
\begin{equation*}
\mathcal{H}\left( V\right) =\Gamma _{hor}\left( OLag,\mathfrak{H}\right) ,
\end{equation*}%
where $\Gamma _{hor}\left( OLag,\mathfrak{H}\right) \subset \Gamma \left(
OLag,\mathfrak{H}\right) $ is the subspace consisting of sections $%
(f_{L^{\circ }}\in \mathcal{H}_{L^{\circ }}:L^{\circ }\in OLag)$ satisfying $%
T_{M^{\circ },L^{\circ }}\left( f_{L^{\circ }}\right) =f_{M^{\circ }}$ for
every $\left( M^{\circ },L^{\circ }\right) \in OLag^{2}$. The vector space $%
\mathcal{H}\left( V\right) $ \ will be referred to as the \textit{canonical
vector space} associated with $V$.

\begin{remark}
The definition of the vector space $\mathcal{H}\left( V\right) $ depends on
a choice of a central character $\psi $.
\end{remark}

The group $Sp(V)$ acts on the set $OLag\times H$ and this induces an action
on the canonical vector space. The representation obtained in this way%
\begin{equation*}
\rho _{V}:Sp(V)\rightarrow GL(\mathcal{H}\left( V\right) ),
\end{equation*}%
will be called the\textit{\ canonical model} of the Weil representation.

\subsection{Construction using kernels of intertwiners}

\subsubsection{Kernel presentation of an intertwining operator}

Fix a pair $\left( M^{\circ },L^{\circ }\right) \in OLag^{2}$ of oriented
Lagrangians and let $%
\mathbb{C}
\left( M\backslash H/L,\psi \right) \subset 
\mathbb{C}
\left( H,\psi \right) $ be the subspace of functions $K\in 
\mathbb{C}
\left( H,\psi \right) $,$\ $satisfying the equivariance property $K\left(
m\cdot h\cdot l\right) =K\left( h\right) $, for every $m\in M$ and $l\in L$.
\ We have an isomorphism of vector spaces 
\begin{equation*}
I:%
\mathbb{C}
\left( M\backslash H/L,\psi \right) \longrightarrow \mathrm{Hom}_{H}(%
\mathcal{H}_{L^{\circ }},\mathcal{H}_{M^{\circ }}),
\end{equation*}%
associating to a function $K\in 
\mathbb{C}
\left( M\backslash H/L,\psi \right) $ the intertwining operator $I\left[ K%
\right] \in $ $\mathrm{Hom}_{H}(\mathcal{H}_{L^{\circ }},\mathcal{H}%
_{M^{\circ }})$ defined by 
\begin{equation*}
I\left[ K\right] \left( f\right) =K\ast f=m_{!}\left( K\boxtimes _{Z\cdot
L}f\right) ,
\end{equation*}%
for\ every $f\in \mathcal{H}_{L^{\circ }}$. Here, $K\boxtimes _{Z\cdot L}f$ $%
\ $denotes the descent of the function $K\boxtimes f\in 
\mathbb{C}
\left( H\times H\right) $ to $H\times _{Z\cdot L}H$---the quotient of $%
H\times H$ by the action $x\cdot (h_{1},h_{2})=(h_{1}x,x^{-1}h_{2})$ for $%
x\in Z\cdot L$---and $m_{!}$ denotes the operation of summation along the
fibers of the multiplication mapping $m:H\times H\rightarrow H$. \ We call
the function $K$ an\textit{\ intertwining kernel}.

It is easy to verify that for a triple $\left( N^{\circ },M^{\circ
},L^{\circ }\right) \in OLag^{3}$ and kernels $K_{2}\in 
\mathbb{C}
\left( N\backslash H/M,\psi \right) $ and $K_{1}\in 
\mathbb{C}
\left( M\backslash H/L,\psi \right) $, their convolution $K_{2}\ast
K_{1}=m_{!}\left( K_{2}\boxtimes _{Z\cdot M}K_{1}\right) $ lies in $%
\mathbb{C}
\left( N\backslash H/L,\psi \right) $. Moreover, the transform $I$ sends
convolution of kernels to composition of operators \ 
\begin{equation*}
I\left[ K_{1}\ast K_{2}\right] =I\left[ K_{1}\right] \circ I\left[ K_{2}%
\right] .
\end{equation*}

\subsubsection{System of canonical intertwining kernels \label%
{sys_kernels_subsub}}

For every $\left( M^{\circ },L^{\circ }\right) \in OLag^{2}$, there exists a
unique kernel $K_{M^{\circ },L^{\circ }}\in $ $%
\mathbb{C}
\left( M\backslash H/L,\psi \right) $ such that $T_{M^{\circ },L^{\circ }}=I%
\left[ K_{M^{\circ },L^{\circ }}\right] $. We will refer to $\left\{
K_{M^{\circ },L^{\circ }}\right\} $ as the system of \textit{canonical
intertwining kernels}. A reformulation of Theorem \ref{SS-vN_thm} is that
the system $\left\{ K_{M^{\circ },L^{\circ }}\right\} $ is multiplicative,
in the sense that $K_{N^{\circ },L^{\circ }}=K_{N^{\circ },M^{\circ }}\ast
K_{M^{\circ },L^{\circ }}$ for every triple $\left( N^{\circ },M^{\circ
},L^{\circ }\right) \in OLag^{3}.$

The system of kernels $\{K_{M^{\circ },L^{\circ }}\}$ can be equivalently
thought of as a single function%
\begin{equation}
K\in 
\mathbb{C}
(OLag^{2}\times H),  \label{K}
\end{equation}
defined by $K\left( M^{\circ },L^{\circ },-\right) =K_{M^{\circ },L^{\circ
}}\left( -\right) $ satisfying\ the following multiplicativity relation on $%
OLag^{3}\times H$ 
\begin{equation}
p_{32}^{\ast }K\ast p_{21}^{\ast }K=p_{31}^{\ast }K\text{,}  \label{mult_eq}
\end{equation}%
where $p_{ji}\left( L_{3}^{\circ },L_{2}^{\circ },L_{1}^{\circ },h\right)
=\left( L_{j}^{\circ },L_{i}^{\circ },h\right) $ are the projections on the $%
j,i$ copies of $OLag$ and the left-hand side of (\ref{mult_eq}) means
fiberwise convolution $p_{32}^{\ast }K\ast p_{21}^{\ast }K(L_{3}^{\circ
},L_{2}^{\circ },L_{1}^{\circ },-)=K\left( L_{3}^{\circ },L_{2}^{\circ
},-\right) \ast K\left( L_{2}^{\circ },L_{1}^{\circ },-\right) $.

\subsubsection{Formula}

In case $M^{\circ },L^{\circ }$ are in generic position, i.e., $\left(
M^{\circ },L^{\circ }\right) \in U_{2}$, the kernel $K_{M^{\circ },L^{\circ
}}$ is given by the following explicit formula%
\begin{equation*}
K_{M^{\circ },L^{\circ }}=A_{M^{\circ },L^{\circ }}\cdot \widetilde{K}%
_{M^{\circ },L^{\circ }},
\end{equation*}%
with $\widetilde{K}_{M^{\circ },L^{\circ }}=\tau ^{\ast }\psi ,$ where $\tau
=\tau _{M^{\circ },L^{\circ }}$ is the inverse of the isomorphism given by
the composition $Z\hookrightarrow H\twoheadrightarrow M\backslash H/L$, and $%
A_{M^{\circ },L^{\circ }}$ as in (\ref{A}).

\subsubsection{Definition of the canonical vector space using kernels\label%
{CVSviaK}}

The canonical vector space $\mathcal{H}(V)$ can be defined using the kernel
function $K$ (\ref{K}). It is the subspace of functions $f\in 
\mathbb{C}
(OLag\times H)$ such that

\begin{itemize}
\item $a_{Z}^{\ast }f=\psi \cdot f,$ where $a_{Z}:Z\times OLag\times
H\rightarrow OLag\times H$ is the action map induced from the action of $Z$
on $H.$

\item $a_{S}^{\ast }f=p^{\ast }f,$ where $S\rightarrow OLag$ is the
tautological vector bundle with fiber $S_{L^{\circ }}=L,$ and $%
p,a_{S}:S\times H\rightarrow OLag\times H$ are the projection and the action
map $a_{S}(l,L^{\circ },h)=(L^{\circ },l\cdot h)$, respectively.

\item We have 
\begin{equation*}
K\ast p_{1}^{\ast }f=p_{2}^{\ast }f,
\end{equation*}%
where $p_{i}:OLag^{2}\times H\rightarrow OLag\times H$ are the projections $%
p_{i}(L_{2}^{\circ },L_{1}^{\circ },h)=(L_{i}^{\circ },h)$.
\end{itemize}

\subsubsection{The Weil representation using the language of kernels\label%
{WRK}}

We consider the set $X=OLag\times H\mathbf{,}$ the action map $%
a_{Sp}:Sp\times X\rightarrow X$, and the pullback operator 
\begin{equation}
a_{Sp}^{\ast }:%
\mathbb{C}
(X)\rightarrow 
\mathbb{C}
(Sp\times X).  \label{ActO}
\end{equation}%
For what follows, it will be convenient for us to denote the canonical
vector space $\mathcal{H(}V\mathcal{)}$ by $\mathcal{H}(X).$ We have the
vector spaces $\mathcal{H}(Sp\times X)$ and $\mathcal{H}(Sp\times Sp\times
X) $ which are defined exactly in the same way (see Subsection \ref{CVSviaK}%
) as $\mathcal{H}(X)$, i.e., by considering the conditions with respect to
the $X$\textbf{-}coordinate\textbf{.} The function\ kernel $K$ satisfies $%
a_{Sp}^{\ast }K=p^{\ast }K$ where $a_{Sp},p:Sp\times OLag^{2}\times
H\rightarrow OLag^{2}\times H$ are the natural action map and projection
map, respectively. Hence, (\ref{ActO}) induces an operator 
\begin{equation*}
\rho _{Sp}:\mathcal{H}(X)\rightarrow \mathcal{H}(Sp\times X).
\end{equation*}%
Note that the homomorphism condition is now manifested via the equality 
\begin{equation*}
(m\times id)^{\ast }\circ \rho _{Sp}=(id\times a_{Sp})^{\ast }\circ \rho
_{Sp},
\end{equation*}%
between the two compositions of the following diagram:%
\begin{equation*}
\begin{diagram} \mathcal{H}(X) & \rTo^{\rho_{Sp}} & \mathcal{H}(Sp \times X)
\\ \dTo^{\rho_{Sp}} & &\dTo_{(id \times a_{Sp})^*} \\ \mathcal{H}( Sp \times
X) & \rTo^{(m \times id)^*} & \mathcal{H}(Sp \times Sp \times X)
\end{diagram}
\end{equation*}%
where $m\mathbf{:}Sp\times Sp\rightarrow Sp$ denotes the multiplication map.

The triple $(\rho _{Sp},Sp,\mathcal{H}(X))$ is trivially identified with the
canonical model of the Weil representation.

\section{Geometric canonical intertwining kernels\label{GCIK}}

In this section, we recall the definition \cite{GH1} of the geometric
counterpart to the set-theoretic system of canonical intertwining kernels.

\subsection{Preliminaries from algebraic geometry}

First, we need to use some space to recall notions and notations from
algebraic geometry and the theory of $\ell $-adic sheaves. \ 

\subsubsection{Varieties}

In the sequel, we are going to translate back and forth between algebraic
varieties defined over the finite field $k$ and their corresponding sets of
rational points. In order to prevent confusion between the two, we use
bold-face letters to denote a variety $\mathbf{X}$ and normal letters $X$ to
denote its corresponding set of rational points $X=\mathbf{X}(k)$. For us, a
variety $\mathbf{X}$ over the finite field is a quasi-projective algebraic
variety, such that the defining equations are given by homogeneous
polynomials with coefficients in the finite field $k$. In this situation,
there exists a (geometric) \textit{Frobenius} endomorphism $Fr:\mathbf{%
X\rightarrow X}$, which is a morphism of algebraic varieties. We denote by $%
X $ \ the set of points fixed by $Fr$, i.e., 
\begin{equation*}
X=\mathbf{X}(k)=\mathbf{X}^{Fr}=\{x\in \mathbf{X}:Fr(x)=x\}.
\end{equation*}

The category of algebraic varieties over $k$ will be denoted by $\mathsf{Var}%
_{k}$.

\subsubsection{Sheaves}

Let $\mathsf{D}^{b}(\mathbf{X)}$ denotes the bounded derived category of
constructible $\ell $-adic sheaves on $\mathbf{X}$ \cite{BBD}. We denote by $%
\mathsf{Perv}(\mathbf{X)}$ the Abelian category of perverse sheaves on the
variety $\mathbf{X}$, that is the heart with respect to the autodual
perverse t-structure in $\mathsf{D}^{b}(\mathbf{X})$. An object $\mathcal{%
F\in }\mathsf{D}^{b}(\mathbf{X)}$ is called $N$-perverse if $\mathcal{F[}%
N]\in \mathsf{Perv}(\mathbf{X)}$. Finally, we recall the notion of a Weil
structure (Frobenius structure) \cite{D3}. A Weil structure associated to an
object $\mathcal{F\in }\mathsf{D}^{b}(\mathbf{X)}$ is an isomorphism%
\begin{equation*}
\vartheta :Fr^{\ast }\mathcal{F}\overset{\sim }{\longrightarrow }\mathcal{F}%
\text{.}
\end{equation*}

A pair $(\mathcal{F},\vartheta )$ is called a Weil object. By an abuse of
notation we often denote $\vartheta $ also by $Fr$. We choose once an
identification $\overline{%
\mathbb{Q}
}_{\ell }\simeq 
\mathbb{C}
$, hence all sheaves are considered over the complex numbers.

\begin{remark}
All the results in this section make perfect sense over the field $\overline{%
\mathbb{Q}
}_{\ell }$, in this respect the identification of $\overline{%
\mathbb{Q}
}_{\ell }$ with $%
\mathbb{C}
$ \ is redundant. The reason it is specified is in order to relate our
results with the standard constructions.
\end{remark}

Given a Weil object $(\mathcal{F},Fr^{\ast }\mathcal{F\simeq F})$ one can
associate to it a function $f^{\mathcal{F}}:X\rightarrow 
\mathbb{C}
$ to $\mathcal{F}$ as follows 
\begin{equation}
f^{\mathcal{F}}(x)=\tsum\limits_{i}(-1)^{i}Tr(Fr_{|H^{i}(\mathcal{F}_{x})}).
\label{S-F}
\end{equation}%
This procedure is called \textit{Grothendieck's sheaf-to-function
correspondence \cite{Gr}}.

\subsection{Geometrization}

We shall now start the geometrization procedure.

\subsubsection{Replacing sets by varieties}

The first step we take is to replace all sets involved by their geometric
counterparts, i.e., algebraic varieties. We denote by $\mathbf{k}$ an
algebraic closure of the field $k$. The symplectic space $(V,\omega )$ is
naturally identified as the set $V=\mathbf{V}(k)$, where $\mathbf{V}$ is a $%
2n$-dimensional symplectic vector space in $\mathsf{Var}_{k}$; the
Heisenberg group $H$ is naturally identified as the set $H=\mathbf{H}(k)$,
where $\mathbf{H}=\mathbf{V\times }\mathbb{G}_{a}$\ is the corresponding
group variety; finally, $OLag=\mathbf{OLag}\left( k\right) $, where $\mathbf{%
OLag}$ is the variety of oriented Lagrangians in $\mathbf{V}$.

\subsubsection{Replacing functions by sheaves}

The second step is to replace functions by their sheaf-theoretic
counterparts \cite{Ga}. The additive character $\psi :k\longrightarrow 
\mathbb{C}
^{\ast }$ is associated via the sheaf-to-function correspondence to the
Artin-Schreier sheaf $\mathcal{L}_{\psi }$ on the variety $\mathbb{G}_{a}$,
i.e., we have $f^{\mathcal{L}_{\psi }}=\psi .$ The Legendre character $%
\sigma $ on $k^{\ast }\simeq $ $\mathbb{G}_{m}(k)$ is associated to the
Kummer sheaf $\mathcal{L}_{\sigma }$ on the variety $\mathbb{G}_{m}$. The
one-dimensional Gauss sum $G(\psi )$ is associated with the Weil object 
\begin{equation*}
\mathcal{G=}\tint\limits_{\mathbb{G}_{a}}\mathcal{L}_{\psi (\frac{1}{2}%
z^{2})}\in \mathsf{D}^{b}(\mathbf{pt}),
\end{equation*}%
where, for the rest of this paper, $\int =\int_{!}$ denotes integration with
compact support \cite{BBD}. Grothendieck's Lefschetz trace formula \cite{Gr}
implies that, indeed, $f^{\mathcal{G}}=G(\psi ).$ In fact, there exists a
quasi-isomorphism $\mathcal{G}\overset{q.i}{\longrightarrow }$ $H^{1}(%
\mathcal{G)}[-1]$ and $\dim H^{1}(\mathcal{G)=}1$, hence, $\mathcal{G}$ can
be thought of as a one-dimensional vector space, equipped with a Frobenius
operator, sitting at cohomological degree $1.$

\subsubsection{Replacing the canonical intertwining kernels by a sheaf}

We recall the geometrization of the function $K$ (\ref{K}) obtained in \cite%
{GH1}.

Let $\mathbf{U}_{2}\mathbf{\subset OLag}^{2}$ be the open subvariety
consisting of pairs $\left( M^{\circ },L^{\circ }\right) \in \mathbf{OLag}%
^{2}$ which are in general position. We define a sheaf \ "of kernels" \ $%
\mathcal{K}_{\mathbf{U}_{2}}$ on the variety $\mathbf{U}_{2}\times \mathbf{H}
$ as follows:

\begin{equation*}
\mathcal{K}_{\mathbf{U}_{2}}=\mathcal{A\otimes }\widetilde{\mathcal{K}}_{%
\mathbf{U}_{2}},
\end{equation*}

where

\begin{itemize}
\item $\widetilde{\mathcal{K}}_{\mathbf{U}_{2}}$ is the sheaf of
non-normalized kernels\textit{\ }given by\textit{\ }%
\begin{equation*}
\widetilde{\mathcal{K}}_{\mathbf{U}_{2}}\left( M^{\circ },L^{\circ }\right) =%
\mathcal{\tau }^{\ast }\mathcal{L}_{\psi },
\end{equation*}%
where $\tau =\tau _{M^{\circ },L^{\circ }}$ is the inverse of the
isomorphism given by the composition 
\begin{equation*}
\mathbf{Z\hookrightarrow H\twoheadrightarrow M\backslash H/L}.
\end{equation*}

\item $\mathcal{A}$ is the "Normalization coefficient" sheaf given by%
\begin{equation*}
\mathcal{A}\left( M^{\circ },L^{\circ }\right) \mathcal{=G}^{\otimes
n}\otimes \mathcal{L}_{\sigma }\left( \left( -1\right) ^{\left( \QATOP{n}{2}%
\right) }\omega _{\wedge }\left( o_{L},o_{M}\right) \right) [2n]\left(
n\right) .
\end{equation*}
\end{itemize}

Let $n_{k}=\dim (\mathbf{OLag}^{k})+n+1$ for $k\in 
\mathbb{N}
$. Consider the projection morphisms $p_{ji}:\mathbf{OLag}^{3}\times \mathbf{%
H\rightarrow }$\ $\mathbf{OLag}^{2}\times \mathbf{H}$ on the $j,i$
coordinates of $\mathbf{OLag}^{3}.$ The main geometric statement of \cite%
{GH1} is the following:

\begin{theorem}[Geometric canonical intertwining kernels]
\label{Gcik}There exist a geometrically irreducible $[n_{2}\mathbf{]}$%
-perverse Weil sheaf $\mathcal{K}$ on $\mathbf{OLag}^{2}\mathbf{\times H}$
of pure weight $w(\mathcal{K})=0$ that satisfies the following three
properties:
\end{theorem}

\begin{enumerate}
\item \textbf{Convolution.}\textit{\ There exists a canonical isomorphism }$%
\theta :p_{32}^{\ast }\mathcal{K}\ast p_{21}^{\ast }\mathcal{K}\simeq
p_{31}^{\ast }\mathcal{K}.$

\item \textbf{Function. }\textit{Applying sheaf-to-function procedure we
recover} $f^{\mathcal{K}}=K.$

\item \textbf{Formula. }\textit{Restricting the sheaf }$\mathcal{K}$\textit{%
\ to the open subvariety} $\mathbf{U}_{2}\mathbf{\times H}$ we have $%
\mathcal{K}_{|\mathbf{U}_{2}\mathbf{\times H}}=\mathcal{K}_{\mathbf{U}_{2}}.$
\end{enumerate}

\section{The idempotent theorem\label{TIT}}

In order to formulate precisely the sign problem and to suggest an answer,
we choose to work in the more general setting of idempotents in categories
with operation \cite{DMOS}.

\subsection{Categorical statement}

Let $\mathcal{C}$ be be a category and \ 
\begin{equation*}
\otimes :\mathcal{C\times C\rightarrow C},\text{ \ \ }(A,B)\mapsto A\otimes
B,\text{ }
\end{equation*}%
a functor. In this paper we will refer to $\otimes $ as an \textit{operation}
in $\mathcal{C}.$

We have two induced functors%
\begin{equation*}
\otimes _{l},\otimes _{r}:\mathcal{C\times C\times C\rightarrow C},
\end{equation*}%
given by $\otimes _{l}(A,B,C)=(A\otimes B)\otimes C$ and $\otimes
_{r}(A,B,C)=A\otimes (B\otimes C).$

An \underline{associativity constraint} for $\otimes $ is a natural
isomorphism 
\begin{equation*}
\alpha :\otimes _{l}\widetilde{\rightarrow }\otimes _{r},
\end{equation*}%
such that, for every $A,B,C,D\in \mathcal{C}$ the diagram

\begin{equation*}
\begin{diagram} ((A \otimes B) \otimes C) \otimes D & \rTo^\alpha & (A
\otimes B) \otimes (C \otimes D) & \rTo^\alpha & A \otimes (B \otimes (C
\otimes D)) \\ \dTo^{\alpha \, \otimes \, id}& & & & \uTo_{id \, \otimes \,
\alpha}\\ (A \otimes (B \otimes C)) \otimes D& & \rTo^{\alpha} & & A \otimes
((B \otimes C) \otimes D) \end{diagram}
\end{equation*}
is commutative. This is usually called the \textit{pentagon coherence }for $%
\alpha $\textit{. }

We would like to study a specific class of objects in $\mathcal{C}$. An 
\underline{idempotent} in $\mathcal{C}$ is a pair $(P,\theta )$ where $P$ is
an object of $\mathcal{C}$ and $\theta $ is an isomorphism 
\begin{equation*}
\theta :P\otimes P\widetilde{\rightarrow }P.
\end{equation*}%
We will also refer to $\theta $ as an idempotent structure on $P.$ We would
like now to study the relation between the idempotent structure and the
associativity constraint.

We have the commutativity

\begin{equation}
\begin{CD} (P\otimes P) \otimes P @>{\alpha}>> P \otimes (P \otimes P)\\
@VV{\theta \, \otimes \, id} V @VV{id \, \otimes \, \theta}V\\ P\otimes P@.
P \otimes P\\ @VV{\theta}V @VV{\theta}V\\ P @>{C}>> P \end{CD}  \label{CDC}
\end{equation}%
where $C:P\rightarrow P$ is by definition the isomorphism which makes the
diagram commutative.

\begin{definition}
We say that the idempotent structure $\theta $ on $P$ is \underline{%
compatible} with the associativity constraint $\alpha $, if in (\ref{CDC})
we have $C=id.$
\end{definition}

The main technical result of the paper is the following:

\begin{theorem}[The idempotent theorem]
\label{ITG}Let $\mathcal{C}$ be a category with an operation $\otimes $
equipped with an associativity constraint $\alpha .$ Let $(P,\theta )$ be an
idempotent in $\mathcal{C}$. Let $C$ be the automorphism of $P$ defined in (%
\ref{CDC}). If $\theta $ satisfies the commutativity of the following
diagrams 
\begin{equation}
\begin{CD} P\otimes P @>{C\otimes id}>> P \otimes P \\ @VV{\theta} V
@VV{\theta}V\\ P @>{C}>> P \end{CD}\text{, \ }\begin{CD} P\otimes P
@>{id\otimes C}>> P \otimes P \\ @VV{\theta} V @VV{\theta}V\\ P @>{C}>> P
\end{CD},  \label{EC}
\end{equation}%
then $\theta $ is compatible with $\alpha $.
\end{theorem}

For a proof of Theorem \ref{ITG} see Appendix \ref{PITG}.\smallskip

For our purposes we will need a straightforward generalization of Theorem %
\ref{ITG}. We consider the sets $\mathbf{2}=\{1,2\},$ $\mathbf{3=}\{1,2,3\},$
$\mathbf{4}=\{1,2,3,4\},$ $\mathbf{5}=\{1,2,3,4,5\},$ and the following
datum:

\begin{itemize}
\item For every $\mathbf{m=2,...,5}$ a category $\mathcal{C}_{\mathbf{m}}$
with an operation $\otimes _{\mathbf{m}}$ equipped with an associativity
constraint $\alpha _{\mathbf{m}}.$

\item For every order preserving embedding $\lambda :\mathbf{l}\rightarrow 
\mathbf{m,}$ $\mathbf{l\varsubsetneq m=3,4,5}$ a faithful functor $%
F_{\lambda }:\mathcal{C}_{\mathbf{l}}\rightarrow \mathcal{C}_{\mathbf{m}}.$
\end{itemize}

We assume that the above functors respect the operations $\otimes _{\mathbf{m%
}}$ and the associativity constraints $\alpha _{\mathbf{m}},$ and are
compatible in the strong sense that the diagram

\begin{equation}
\begin{diagram} \mathcal{C}_{\bf 3} &\rTo^{F_\sigma} & \mathcal{C}_{\bf 4} &
\rTo^{F_\tau}& \mathcal{C}_{\bf 5}\\ &\luTo_{F_\lambda} &\uTo_{F_\mu}
&\ruTo_{F_\nu} \\ & &\mathcal{C}_{\bf 2} \end{diagram}  \label{Comp}
\end{equation}%
is commutative for $\mu =\sigma \circ \lambda $ and $\upsilon =\tau \circ
\mu .$

For simplicity, let us denote the above datum by $(\mathcal{C}_{\bullet
},\otimes _{\bullet },\alpha _{\bullet },F_{\bullet }).$ Then we can define
the notion of an $F_{\bullet }$-idempotent. This is an object $P\in \mathcal{%
C}_{\mathbf{2}}$ equipped with an isomorphism in $\mathcal{C}_{\mathbf{3}}$ 
\begin{equation*}
\theta :F_{32}(P)\otimes _{\mathbf{3}}F_{21}(P)\widetilde{\rightarrow }%
F_{31}(P),
\end{equation*}%
where by $ji$ we mean the function from $\mathbf{2}$ to $\mathbf{3}$ given
by $ji(1)=i$ and $ji(2)=j,$ $1\leq i<j\leq 3.$ For such idempotent we have
the commutative diagram in $\mathcal{C}_{\mathbf{4}}$ 
\begin{equation}
\begin{CD} (F_{43}(P)\otimes F_{32}(P)) \otimes F_{21}(P) @>{\alpha}>>
F_{43}(P) \otimes (F_{32}(P) \otimes F_{21}(P))\\ @VV{\theta \, \otimes \,
id} V @VV{id \, \otimes \, \theta}V\\ F_{42}(P)\otimes F_{21}(P)@. F_{43}(P)
\otimes F_{31}(P)\\ @VV{\theta}V @VV{\theta}V\\ F_{41}(P) @>{C}>> F_{41}(P)
\end{CD}  \label{IT-G}
\end{equation}%
where $ji$ means the function from $\mathbf{2}$ to $\mathbf{4}$ given by $%
ji(1)=i$ and $ji(2)=j,$ $1\leq i<j\leq 4,$ $\otimes =\otimes _{\mathbf{4},}$ 
$\alpha =\alpha _{\mathbf{4}},$ and $\theta $ is an abbreviation for the
various isomorphisms $F_{\sigma }(\theta )$ for the appropriate $\sigma $'s.

Exactly the same argument as in the proof of Theorem \ref{ITG} yields the
following theorem:

\begin{theorem}[The idempotent theorem---generalization]
\label{Tit-G}Let $(\mathcal{C}_{\bullet },\otimes _{\bullet },\alpha
_{\bullet },F_{\bullet })$ be a datum as defined above. Let $(P,\theta )$ be
an $F_{\bullet }$-idempotent. If $\theta $ satisfies the analogue of (\ref%
{EC}), then $\theta $ is compatible with $\alpha _{\bullet },$ i.e., in (\ref%
{IT-G}) we have $C=id.$
\end{theorem}

\section{Applications\label{Appl}}

We are ready to formulate and solve the sign problem, and to suggest several
additional applications. The main application is the definition of the
canonical category, which yields the canonical model of the categorical Weil
representation.

\subsection{Application I: Solution to the sign problem}

We follow the setup of Sections \ref{GCIK} and \ref{TIT}. Consider the
following datum:

\begin{itemize}
\item For every $\mathbf{m=2,...,5}$ the category $\mathcal{C}_{\mathbf{m}}=%
\mathsf{D}^{b}(\mathbf{OLag}^{m}\mathbf{\times H})$ with the operation $%
\otimes _{\mathbf{m}}$which is induced by the convolution $\ast $ on $%
\mathbf{H,}$ and with the standard associativity constraint $\alpha _{%
\mathbf{m}}.$

\item For every order preserving embedding $\lambda :\mathbf{l}\rightarrow 
\mathbf{m,}$ $\mathbf{l\varsubsetneq m=3,4,5}$, a faithful functor $%
F_{\lambda }:\mathcal{C}_{\mathbf{l}}\rightarrow \mathcal{C}_{\mathbf{m}}$
given by the pullback functor $p_{\lambda }^{\ast },$ where $p_{\lambda }:$ $%
\mathbf{OLag}^{m}\mathbf{\times H\rightarrow OLag}^{l}\mathbf{\times H}$ is
the projection morphism induced by $\lambda .$
\end{itemize}

The datum $(\mathcal{C}_{\bullet },\otimes _{\bullet },\alpha _{\bullet
},F_{\bullet })$ satisfies the compatibility of Diagram (\ref{Comp}). In
these terms the pair $(\mathcal{K},\theta \mathcal{)}$ where $\mathcal{K}$
is the sheaf of canonical geometric intertwining kernels (see Theorem \ref%
{Gcik}) is an $F_{\mathbf{\bullet }}$-idempotent. Recall that $\mathcal{K}$
is geometrically irreducible. The precise formulation of the sign problem 
\cite{B, D2} is the following:

\begin{problem}[The sign problem]
Compute the value of the scalar $c$ in the morphism $C=c\cdot id$ that
appears in diagram (\ref{IT-G}) with $P$ replaced by the sheaf $\mathcal{K}.$
\end{problem}

Theorem \ref{Tit-G} suggests the answer.

\begin{theorem}[The idempotent theorem---particular case]
\label{Titpc}We have $c=1.$
\end{theorem}

\subsection{Application II: The canonical category\label{TCC}}

Theorem \ref{Titpc} suggests the definition of a canonical category $%
\mathcal{C}(\mathbf{V})$ associated with a symplectic vector space variety $%
\mathbf{(V,\omega )}$ defined over $k.$ This category, and its definition,
are the geometric analogue of the canonical vector space $\mathcal{H}(V),$
and its definition, via the kernel function $K,$ as described in Subsection %
\ref{CVSviaK}.

The canonical category $\mathcal{C}(\mathbf{V})$ consists of sheaves $%
\mathcal{F}\in \mathsf{D}^{b}(\mathbf{OLag\times H})$ equipped with the
following structure:

\begin{itemize}
\item An isomorphism $\iota _{\mathbf{Z}}:a_{\mathbf{Z}}^{\ast }\mathcal{F}%
\widetilde{\rightarrow }\mathcal{L}_{\psi }\boxtimes \mathcal{F},$ where $a_{%
\mathbf{Z}}:\mathbf{Z\times OLag\times H}\rightarrow \mathbf{OLag\times H}$
is the action morphism induced from the action of $\mathbf{Z}$ on $\mathbf{H}%
.$

\item An isomorphism $\iota _{\mathbf{S}}^{\ast }:a_{\mathbf{S}}^{\ast }%
\mathcal{F}\widetilde{\rightarrow }$ $p^{\ast }\mathcal{F},$ where $\mathbf{%
S\rightarrow OLag}$ is the tautological vector bundle with fiber $\mathbf{S}%
_{\mathbf{L}^{\circ }}=\mathbf{L},$ and $p,a_{\mathbf{S}}\mathbf{:S\times
H\rightarrow OLag\times H}$ are the projection morphism and the action
morphism $a_{\mathbf{S}}(l,\mathbf{L}^{\circ },h)=(\mathbf{L}^{\circ
},l\cdot h)$, respectively.

\item An isomorphism 
\begin{equation*}
\eta :\mathcal{K}\ast p_{1}^{\ast }\mathcal{F}\widetilde{\rightarrow }%
p_{2}^{\ast }\mathcal{F},
\end{equation*}%
where $p_{i}:\mathbf{OLag}^{2}\mathbf{\times H\rightarrow OLag\times H}$, $%
i=1,2$, are the projection morphisms $p_{i}(\mathbf{L}_{2}^{\circ },\mathbf{L%
}_{1}^{\circ },h)=(\mathbf{L}_{i}^{\circ },h)$, and $\mathcal{K}$ is the
sheaf of canonical geometric intertwining kernels (Theorem \ref{Gcik}).
\end{itemize}

The isomorphism $\eta $ is required to be compatible with the associativity
constraint $\alpha ,$ induced from the convolution $\ast $ on $\mathbf{H,}$
and the idempotent structure $\theta $ of $\mathcal{K}$, i.e., the following
diagram 
\begin{equation*}
\begin{CD} (p^*_{32}\mathcal{K} \ast p^*_{21}\mathcal{K}) \ast
p^*_1\mathcal{F} @>{\alpha}>> p^*_{32}\mathcal{K} \ast (p^*_{21}\mathcal{K}
\ast p_1\mathcal{F})\\ @VV{\theta \, \ast \, id} V @VV{id \, \ast \,
\eta}V\\ p^*_{31}\mathcal{K} \ast p^*_1\mathcal{F} @. p^*_{32}\mathcal{K}
\ast p^*_2\mathcal{F}\\ @VV{\eta}V@VV{\eta}V\\ p^*_3\mathcal{F}@>{id}>>
p^*_3\mathcal{F} \end{CD}
\end{equation*}%
is required to be commutative. Here, $p_{ji}:\mathbf{OLag}^{3}\mathbf{\times
H\rightarrow OLag}^{2}\mathbf{\times H,}$ $1\leq i<j\leq 3,$ and $p_{i}:%
\mathbf{OLag}^{3}\mathbf{\times H\rightarrow OLag\times H,}$ $1\leq i\leq 3,$
denote the projection morphisms on the $ji$, and $i$ coordinates,
respectively.

\begin{remark}[Important remark]
The category $\mathcal{C}(\mathbf{V})$ makes sense, i.e., it is nontrivial,
only if $c=1$ in Theorem \ref{Titpc}. The verification of this assertion,
follows, more or less the same argument as in the proof of the idempotent
theorem, and hence will be left for the reader.
\end{remark}

\begin{remark}
Let us denote by $\mathbf{Symp}$ the (groupoid) category whose objects are
finite dimensional symplectic vector spaces in the category of algebraic
varieties over $k$, and morphisms are linear isomorphisms of symplectic
vector spaces. It is not hard to see that for each morphism $f:\mathbf{%
U\rightarrow V}$ in $\mathbf{Symp,}$ we have an induced pullback functor $%
f^{\ast }:$\ \ $\mathcal{C}(\mathbf{V})\rightarrow \mathcal{C}(\mathbf{U})$.
Moreover, for any pair of morphisms $\mathbf{U}\overset{f}{\mathbf{%
\rightarrow }}\mathbf{V}\overset{g}{\mathbf{\rightarrow }}\mathbf{W}$ we
have, in this case, the equality ($g\circ f)^{\ast }=f^{\ast }\circ g^{\ast
}.$ This means that we obtain the categorical analog of the functor (\ref{Q}%
), i.e., we constructed an explicit quantization lax $2$-functor \cite{Gr1,
V} 
\begin{equation*}
\mathbf{V\mapsto }\text{ }\mathcal{C}(\mathbf{V}),
\end{equation*}%
associating the category of $\ell $-adic sheaves $\mathcal{C}(\mathbf{V})$
to the object $\mathbf{V}$.
\end{remark}

\subsubsection{The subcategories of perverse and Weil objects\label{SC}}

Let us consider the subcategories $\mathsf{Perv}(\mathbf{OLag\times H)}$ and 
$\mathsf{D}_{\mathsf{w}}^{b}(\mathbf{OLag\times H}),$ of $\mathsf{D}^{b}(%
\mathbf{OLag\times H}),$ consisting of perverse sheaves, and Weil sheaves,
respectively (see Section \ref{GCIK}). These categories are preserved by the
action of the kernel $\mathcal{K}$. For the case of $\mathsf{Perv}(\mathbf{%
OLag\times H)}$, since $\mathcal{K}$ is essentially a Fourier transform
kernel, this follows from the Katz--Laumon theorem \cite{KL}. Hence, it
makes sense to define in $\mathcal{C}(\mathbf{V})$ the subcategories $%
\mathsf{P(}\mathbf{V)}$, $\mathcal{C}_{\mathsf{w}}(\mathbf{V})$, and $%
\mathsf{P}_{\mathsf{w}}\mathsf{(}\mathbf{V)=}$ $\mathsf{P(}\mathbf{V)}\cap 
\mathcal{C}_{\mathsf{w}}(\mathbf{V}),$ of perverse sheaves, Weil sheaves,
and perverse Weil sheaves, respectively, which satisfies the conditions of
Subsection \ref{TCC} above. In particular, using the sheaf-to-function
procedure (\ref{S-F}) we can recover the vector space $\mathcal{H}(V)$ from
the categories $\mathcal{C}_{\mathsf{w}}(\mathbf{V})$, or $\mathsf{P}_{%
\mathsf{w}}\mathsf{(}\mathbf{V).}$ This means that we obtained the
geometrizations of the vector space $\mathcal{H}(V)$ by these categories.

\subsection{Application III: The categorical Weil representation}

The canonical vector space $\mathcal{H}(V)$ yields a natural model for the
Weil representation of $Sp$. In the same spirit, the canonical category $%
\mathcal{C}(\mathbf{V})$ yields a model of, what we will call, the \textit{%
categorical Weil representation} of the algebraic group $\mathbf{Sp=}Sp(%
\mathbf{V,\omega }).$

We consider the variety $\mathbf{X=OLag\times H,}$ the action morphism $a_{%
\mathbf{Sp}}:\mathbf{Sp\times X\rightarrow X}$, and the pullback functor 
\begin{equation}
a_{\mathbf{Sp}}^{\ast }:\mathsf{D}^{b}(\mathbf{X})\rightarrow \mathsf{D}^{b}(%
\mathbf{Sp\times X}).  \label{ActD}
\end{equation}%
For what follows it will be convenient for us to denote the category $%
\mathcal{C}(\mathbf{V})$ by $\mathcal{C}(\mathbf{X}).$ We have the
categories $\mathcal{C}(\mathbf{Sp\times X})$ and $\mathcal{C}(\mathbf{%
Sp\times Sp\times X})$ which are defined exactly in the same way (see
Subsection \ref{TCC}) as $\mathcal{C}(\mathbf{X})$, i.e., by considering the
conditions with respect to the $\mathbf{X}$\textbf{-}coordinate\textbf{.}
The sheaf\ $\mathcal{K}$ of geometric intertwining kernels (Theorem \ref%
{Gcik}) is equipped with a natural isomorphism $\iota _{\mathbf{Sp}}:a_{%
\mathbf{Sp}}^{\ast }\mathcal{K}\widetilde{\rightarrow }p^{\ast }\mathcal{K}$
where $a_{\mathbf{Sp}},p:\mathbf{Sp\times OLag}^{2}\times \mathbf{%
H\rightarrow OLag}^{2}\times \mathbf{H}$ are the natural action morphism and
projection morphism, respectively. Hence, the functor (\ref{ActD}) induces a
functor 
\begin{equation}
\mathbf{\rho }_{\mathbf{Sp}}:\mathcal{C}(\mathbf{X})\rightarrow \mathcal{C}(%
\mathbf{Sp\times X}).  \label{CWR}
\end{equation}%
Note that, in our specific situation, we have an equality 
\begin{equation*}
(m\times id)^{\ast }\circ \mathbf{\rho }_{\mathbf{Sp}}=(id\times a_{\mathbf{%
Sp}})^{\ast }\circ \mathbf{\rho }_{\mathbf{Sp}},
\end{equation*}%
between the two compositions of the following diagram:%
\begin{equation*}
\begin{diagram} \mathcal{C}(\bf X) & \rTo^{\bf{\rho_{Sp}}} & \mathcal{C}(\bf
Sp \times X) \\ \dTo^{\bf{\rho_{Sp}}} & &\dTo_{(id \times a_{\bf Sp})^*} \\
\mathcal{C}(\bf Sp \times X) & \rTo^{(m \times id)^*} & \mathcal{C}(\bf Sp
\times Sp \times X) \end{diagram}
\end{equation*}%
where $m\mathbf{:Sp\times Sp\rightarrow Sp}$ denotes the multiplication
morphism.

We will call the triple $(\mathbf{\rho }_{\mathbf{Sp}},\mathbf{Sp,}\mathcal{C%
}(\mathbf{X}))$ the \textit{canonical model of the categorical Weil
representation. }

\subsubsection{ Action on subcategories}

The categorical Weil representation $(\mathbf{\rho }_{\mathbf{Sp}},\mathbf{%
Sp,}\mathcal{C}(\mathbf{X}))$ induces the categorical representations $(%
\mathbf{\rho }_{\mathbf{Sp}},\mathbf{Sp,}\mathsf{P}(\mathbf{X}))),$ $(%
\mathbf{\rho }_{\mathbf{Sp}},\mathbf{Sp,}\mathcal{C}_{\mathsf{w}}(\mathbf{X}%
)),$ and $(\mathbf{\rho }_{\mathbf{Sp}},\mathbf{Sp,}\mathsf{P}_{\mathsf{w}}(%
\mathbf{X}))$ by restricting the functor $\mathbf{\rho }_{\mathbf{Sp}}$ (\ref%
{CWR}) to the various subcategories introduced in Subsection \ref{SC}. In
particular, using the sheaf-to-function procedure (\ref{S-F}) we can recover
the canonical model $(\rho _{Sp},Sp,\mathcal{H}(X))$ of the Weil
representation (see Subsection \ref{WRK}) from the categorical
representations $(\mathbf{\rho }_{\mathbf{Sp}},\mathbf{Sp,}\mathcal{C}_{%
\mathsf{w}}(\mathbf{X}))$ and $(\mathbf{\rho }_{\mathbf{Sp}},\mathbf{Sp,}%
\mathsf{P}_{\mathsf{w}}(\mathbf{X}))\mathbf{.}$ This means that we obtained
the geometrizations of the representation $(\rho _{Sp},Sp,\mathcal{H}(X))$
by these categorical representations.

\subsection{Application IV: The geometric Weil representation}

In \cite{GH3} the authors elaborated on \cite{D1}, and developed the theory
of invariant presentation of the Weil representation and its geometric
analogue.

Consider the Heisenberg representation $(\pi ,H,\mathcal{H})$ (Subsection %
\ref{HR}), and the associated Weil representation $(\rho ,Sp,\mathcal{H})$.
We define the kernel of the Weil representation as the function $K(g,h)=%
\frac{1}{\dim \mathcal{H}}\cdot Tr(\rho (g)\circ \pi (h^{-1}))$, for every $%
g\in Sp,$ $h\in H.$ We have $K(g,\bullet )\in 
\mathbb{C}
(H,\psi ^{-1})$---the space of complex valued functions $f$ on $H$ such that 
$f(z\cdot h)=\psi ^{-1}(z).$ Hence, we can, and will, consider $K$ as a
function 
\begin{equation*}
K:Sp\times V\rightarrow 
\mathbb{C}
.
\end{equation*}%
The function $K$ determines the Weil representation $\rho $ completely,
since we have $\pi (K(g,\bullet )):=\underset{v\in V}{\tsum }K(g,v)\cdot \pi
(v)=\rho (g),$ and it satisfies the convolution property 
\begin{equation*}
p_{1}^{\ast }K\mathcal{\ast }p_{2}^{\ast }K\mathcal{=(}m\times id)^{\ast }K,
\end{equation*}%
where $m:Sp\times Sp\rightarrow Sp$ denotes the multiplication map, and $%
p_{1},p_{2}:$ $Sp\times Sp\times V\rightarrow Sp\times V$ denote the
projections on $Sp\times V$ via the first and second $Sp$ coordinates,
respectively, and the $\ast $ denotes the operation on the function space $%
\mathbb{C}
(H,\psi ^{-1})=%
\mathbb{C}
(V)$ induced from the convolution operation on functions on the Heisenberg
group $H$.

The two main contributions of\ \cite{GH3} are:

\begin{itemize}
\item An explicit formula for the function $K$ as follows:%
\begin{equation*}
K\left( g,v\right) =\tfrac{1}{\dim \mathcal{H}}\cdot \sigma \left(
(-1)^{n}\cdot \det \left( g-I\right) \right) \cdot \psi \left( \tfrac{1}{4}%
\omega \left( \kappa \left( g\right) v,v\right) \right)
\end{equation*}%
for every $g\in Sp$ such that $g-I$ is invertible, where $n=\dim (V),$ and $%
\sigma $ denotes the unique nontrivial quadratic character (Legendre
character) of the multiplicative group $k^{\ast },$ and $\kappa \left(
g\right) =\frac{g+I}{g-I}$ is the Cayley transform.

\item The construction of the sheaf theoretic object that geometrizes the
function $K$. This is the content of the following theorem (see Section \ref%
{GCIK} for notations and definitions) :
\end{itemize}

\begin{theorem}[Geometric Weil representation]
There exists a geometrically irreducible $[\dim \mathbf{Sp}]$-perverse Weil
sheaf $\mathcal{K}$ of pure weight zero on $\mathbf{Sp\times V}$ satisfying
the following properties:

\begin{enumerate}
\item \textbf{Convolution}.\textbf{\ }There exists a canonical isomorphism $%
\theta :p_{1}^{\ast }\mathcal{K\ast }p_{2}^{\ast }\mathcal{K}\widetilde{%
\mathcal{\rightarrow }}(m\times id)^{\ast }\mathcal{K}.$

\item \textbf{Function. }We have $f^{\mathcal{K}}=K.$

\item \textbf{Formula}. For every $g\mathbf{\in Sp}$ with $\det (g-I\mathbf{%
)\neq }0$ we have 
\begin{equation*}
\mathcal{K(}g,v)=\mathcal{L}_{\sigma }((-1)^{n}\cdot \det (g-I))\otimes 
\mathcal{L}_{\psi }\left( \tfrac{1}{4}\omega (\kappa (g)v,v)\right) [2n](n).
\end{equation*}
\end{enumerate}
\end{theorem}

Here $m:\mathbf{Sp\times Sp\rightarrow Sp}$ denotes the multiplication
morphism, and $p_{1},p_{2}:$ $\mathbf{Sp\times Sp\times V\rightarrow
Sp\times V}$ denote the projections on $\mathbf{Sp\times V}$ via the first
and second $\mathbf{Sp}$ coordinates, respectively.

The following is similar to the sign problem. Consider the commutative
diagram of isomorphisms in $\mathsf{D}^{b}(\mathbf{Sp}^{3}\mathbf{\times V})$

\begin{equation*}
\begin{CD} (p^*_1\mathcal{K} \ast p^*_2\mathcal{K}) \ast p^*_3\mathcal{K}
@>{\alpha}>> p^*_1\mathcal{K} \ast (p^*_2\mathcal{K} \ast
p^*_3\mathcal{K})\\ @VV{\theta \, \ast \, id} V @VV{id \, \ast \, \theta}V\\
(m_{12} \times id)^*\mathcal{K} \ast p^*_3\mathcal{K} @. p^*_1\mathcal{K}
\ast (m_{23} \times id)^*\mathcal{K}\\ @VV{\theta}V @VV{\theta}V\\
m^*_{123}\mathcal{K} @>{C}>> m^*_{123}\mathcal{K} \end{CD}
\end{equation*}%
where $\alpha $ is the \textit{associativity constraint}\ \cite{DMOS} for
convolution $\ast $ of sheaves on the Heisenberg group $\mathbf{H}$, $\theta 
$ is the isomorphism appearing in the convolution property of the sheaf $%
\mathcal{K}$, and $C$ is by definition the isomorphism that makes the
diagram commutative. In addition $m_{ij,}m_{123}:\mathbf{Sp}^{3}\rightarrow 
\mathbf{Sp}$ are the morphisms defined by multiplication of the $i,j$
coordinates, $1\leq i<j\leq 3,$ and all the three coordinates, respectively,
and $p_{1},p_{2}:$ $\mathbf{Sp\times Sp\times V\rightarrow Sp\times V}$
denote the projections on $\mathbf{Sp\times V}$ via the first and second $%
\mathbf{Sp}$ coordinates.

The sheaf $\mathcal{K}$ is geometrically irreducible, hence, $C=c\cdot id$
is a scalar morphism.\smallskip

\begin{problem}[The sign problem]
Compute the scalar $c.$
\end{problem}

A variant of Theorem \ref{Tit-G}, with the same proof, can be applied also
in this case, yielding:

\begin{theorem}[The idempotent theorem---particular variant]
We have $c=1.$
\end{theorem}

\appendix

\section{Proof of Theorem \protect\ref{ITG}\label{PITG}}

We have the commutativity

\begin{equation}
\begin{diagram} ((P \otimes P) \otimes P) \otimes P & \rTo^\alpha & (P
\otimes P) \otimes (P \otimes P) & \rTo^\alpha & P \otimes (P \otimes (P
\otimes P)) \\ \dTo^{\alpha \, \otimes \, id}& & & & \uTo_{id \, \otimes \,
\alpha}\\ (P \otimes (P \otimes P)) \otimes P& & \rTo^{\alpha} & & P \otimes
((P \otimes P) \otimes P) \end{diagram}  \label{CDA}
\end{equation}%
By successive application of $\theta $, each term in (\ref{CDA}) is
identified with $P,$ and using (\ref{CDC}), (\ref{EC}), and the naturality
of $\alpha ,$ the arrows become

\begin{equation*}
\begin{diagram} P & \rTo^C & P & \rTo^C & P \\ \dTo^ C& & & & \uTo_C\\ P& &
\rTo^C & & P \end{diagram}
\end{equation*}%
which implies that $C^{3}=C^{2},$ i.e., $C=id$ as claimed.

\end{document}